\font\eightrm=cmr8  
\font\eighttt=cmtt8
\magnification=\magstephalf
\overfullrule=0in
 
\bf
\noindent
A SHORT PROOF OF JACOBI'S FORMULA FOR THE NUMBER
OF REPRESENTATIONS OF AN INTEGER AS A SUM OF FOUR SQUARES
\medskip
\noindent
\it
George E. Andrews\footnote{$^1$}
{\eightrm 
Dept. of Mathematics, Pennsylvania State University,
University Park, PA 16182, USA. 
\break
{\eighttt andrews@cantor.math.psu.edu .}
Supported in part by the NSF.
},
Shalosh B. Ekhad\footnote{$\,^2$}
{
\eightrm Dept. of Mathematics, Temple University,
Philadelphia, PA 19122, USA. 
\break
{ \eighttt ekhad@euclid.math.temple.edu ,
 zeilberg@euclid.math.temple.edu .}
Supported in part by the NSF.
}
and Doron Zeilberger{$\,{}^2$}
\rm
\medskip
(Appeard in the Amer. Math. Monthly 100(1993), 273-276)
\bigskip
Diophantus probably knew, and Lagrange[L]
proved, that every positive integer can be written as a sum
of four perfect squares. Jacobi[J] proved the stronger
result  that the number
of ways in which a positive integer can be so written\footnote{$^3$}
{ \eightrm Precisely: the number of vectors ({\it not sets})
$( x_1 , x_2 , x_3 , x_4 )$, where
the components are (positive, negative, or zero) integers,
such that 
$ x_1^2 + x_2^2 + x_3^2 + x_4^2 =n$.}
equals $8$ times the sum of its divisors
that are not multiples of 4. Here we give a short new
proof that only uses high school algebra, and is completely
{\it from scratch.}
All infinite series and products
that appear are to be taken in the entirely elementary sense of
formal power series.

The problem of representing integers as sums of squares
has drawn the attention of many great mathematicians, and we encourage
the reader to look up Grosswald's[G] erudite masterpiece on this
subject.

The crucial part of our proof is played by two simple
identities, that we state as one Lemma.
\medskip
\noindent
{\bf Lemma :} Let
\medskip
$$
H_n = H_n (q) = {{1+q} \over {1-q}} {{1+ q^2} \over {1- q^2}}
... {{1+ q^n} \over {1- q^n }} .
\hbox { {\rm For all integers $n \geq 0$.}}
$$
\medskip
$$
\sum_{k=-n}^{n}
{{4 (-q)^k } \over {(1+ q^k )^2 }} H_n^{2}  H_{n+k}
 H_{n-k} =1, 
\eqno(a)
$$
\medskip
$$
\sum_{k=0}^{n} {{2 (- q^{n+1} )^k } \over {1+ q^k}}
{{H_k} \over {H_n}} =
\sum_{k=-n}^{n} (-q)^{k^2} .
\eqno(b)
$$
\medskip
\noindent
{\bf Proof:} Let $L_1 (n)$ and $L_2 (n)$ be the left sides
of (a) and (b) respectively, and let $F_1 (n,k)$, and
$F_2 (n,k)$ be the respective summands. Since both
(a) and (b) obviously hold for $n=0$, it suffices to prove
that for every $n \geq 0$, $L_1 (n+1 ) - L_1 (n) =0$, and 
$L_2 (n+1) - L_2 (n) = 2 (-q)^{(n+1)^2} $.
To this end, we construct
\medskip
$$
G_1 (n,k):={ {q^{n-k+1} (1+ q^{2n+2} ) (1+ q^k )^2
(1+ q^{n+k+1} )} \over
{(1- q^{n+1} )^3 (1- q^{n+k+1} )(1+ q^{n+1} ) }} F_1 (n,k) ,\;
G_2 (n,k) :=
{ {(- q^{n+1}) (1+ q^k )} \over  
{1 + q^{n+1} }}
  F_2 (n,k),
$$
\medskip
\noindent
with the motive that
\smallskip
$$
F_1 (n+1,k) - F_1 (n,k)= G_1 (n,k)- G_1 (n,k-1), \;
F_2 (n+1,k)- F_2 (n,k)= G_2 (n,k)- G_2 (n,k-1) ,
\eqno(1)
$$
\medskip
\noindent
which immediately imply them, by telescoping,  upon
summing from $k=-n-1$ to $k=n+1$,
and from $k=0$ to $k=n+1$ respectively.
The two identities of (1) are purely routine, since
dividing through by $F_1 (n,k)$ and $F_2 (n,k)$ respectively,
lead to routinely-verifiable high-school-algebra identities.\quad

Dividing both sides of (a) by $H_n^4$ and letting
$n \rightarrow \infty$ in (a) and (b) gives
\medskip
$$
1+ 8 \sum_{k=1}^{\infty}{ { (-q)^k } \over {(1+ q^k )^2 }}
= H_{\infty}^{-4} ,
\eqno(a')
$$
\medskip
$$
\sum_{k= - \infty }^{\infty} (-q)^{k^2} =
H_{\infty}^{-1} .
\eqno(b')
$$
\medskip
\noindent
Combining (a') and (b'), yields, after
changing $q \rightarrow -q$,
\medskip
$$
( \sum_{k= - \infty }^{\infty} q^{ k^2 } \, )^4 =
1+ 8  \sum_{k=1}^{\infty} {{q^k} \over { (1 + (-q)^k )^2 }} \; .
\eqno(2)
$$
\medskip
\noindent
The coefficient of a typical term $q^n$ on the left of (2)
is the number of ways of writing $n$ as a sum of
four squares. It remains to show that the coefficient of $q^n$
in the sum on the right of (2) equals the sum of the divisors of $n$
that are not multiples of $4$.
 
Using the power-series expansion
$z/ (1+z)^2 =  \sum_{r=1}^{\infty} (-1)^{(r+1)} r z^r $, with
$z= (-q)^k$, and collecting like powers,
the sum on the right side may be rewritten
\medskip
$$
\sum_{k=1}^{\infty} \sum_{r=1}^{\infty} (-1)^{(k+1)(r+1)} \, r \,  q^{kr} =
\sum_{n=1}^{\infty} q^n \, [ \sum_{r \mid n}  (-1)^{(r+1)(n/r+1)} r \, ].
$$
\medskip
\noindent
The coefficient of $q^n$ above is a weighed sum of divisors
$r$ of $n$, where the coefficient of $r$ is $-1 = +1 -2$ if both $r$ and
$n/r$ are even and $+1$ otherwise, so the coefficient of $q^n$ is
\medskip
$$
\sum_{r \mid n} r - \sum_{  {{ r \mid n} \atop
{ r , n/r {\, even}} }} 2r
=
\sum_{ d \mid n} d - \sum_{ {{ d \mid n} \atop {4 \mid d} }} d
=
 \sum_{ {{ d \mid n} \atop {4 \not{\,\mid} d } }} d \, .
$$
\bigskip
The finitary identities (a) and (b) combine to yield
a single finitary identity
\medskip
$$
\left ( {\sum_{k=0}^{n} {{2 (- q^{n+1} )^k } \over {1+ q^k}}
{H_k} } \right )^4 \,
\sum_{k=-n}^n 
{{4 (-q)^k } \over {(1+ q^k )^2 }} { { H_{n+k} } \over { H_n} }
{{ H_{n-k} } \over { H_n }}
=  ( \sum_{k= -n }^n (-q)^{ k^2 } )^4 \; ,
\eqno(3)
$$
\medskip
\noindent
which also immediately implies Jacobi's theorem, by taking it 
``mod $q^n $''  for any desired $n$. Identity (3) makes
it transparent that our proof only uses the potential infinity,
not the ultimate one.

The identities of the Lemma  are examples of q-binomial coefficient 
identities, a.k.a terminating basic hypergeometric series identities. The
proof of such identities is now completely routine[WZ][Z]. The
proof of the Lemma given here used the algorithm of [Z]. Further
applications of basic hypergeometric series to number theory can be
found in [A1]. An excellent modern reference to basic hypergeometric 
series is [GR].

{\eightrm
We conclude with some comments addressed mainly to the 
cognoscenti.
Identities (a) and (b) are special cases of classical
identities: (a) is a special case of Jackson's theorem 
[GR, p. 35, eq. (2.6.2)], and
(b) is a special case of Watson's q-analog of Whipple's theorem 
([GR, p.35, eq. (2.5.1)], see also [A2, p. 118, eq. (4.3)]. )
The discovery of (b) was motivated by [S1] and [S2].
 
We see fairly clearly how
to do the 2-square theorem (a different instance of Jackson's theorem
replaces (a)); however the theorems for 6 and 8 squares apparently require
(using this approach) some instance of the  
$ _6 \Psi_6 $ summation theorem [GR,
p. 128, (5.3.1)] (see [A1, pp. 461-465] for details).  
Since we do not know
a finitary analog of the $_6 \Psi_6 $ summation, 
the question of a similar result for 6 and 8 squares is of interest.}
\medskip
\noindent
{\bf Acknowledgement:} The referee made several helpful comments.
\bigskip
\noindent
{\bf References}
\medskip
\noindent
[A1] G.E. Andrews, {\it Applications of basic hypergeometric functions},
SIAM Rev. {\bf 16}(1974), 441-484.
 
\noindent
[A2] G.E. Andrews, {\it The fifth and seventh order mock theta 
functions}, Trans. Amer. Math. Soc., {\bf 293} (1986), 113-134.
 
\noindent
[GR] G. Gasper and M. Rahman, "{\it Basic hypergeometric series}",
Cambridge University Press, 1990.
 
\noindent
[G] E. Grosswald, "{\it Representations of 
integers as sums of squares}", Springer, New York, 1985.
 
\noindent
[J] C.G.J. Jacobi, {\it Note sur la d\'ecomposition d'un nombre donn\'e
en quatre carre\'es} , J. Reine Angew. Math. {\bf 3} (1828), 191.
{\it Werke}, vol. I, 247.
 
\noindent
[L] J.L. Lagrange, Nouveau M\'em. Acad. Roy. Sci. Berlin(1772),
123-133; {\it Oevres}, vol. 3, 189-201.
 
\noindent
[S1] D. Shanks, {\it A short proof of an identity of Euler}, 
Proc. Amer. Math. Soc., {\bf 2} (1951),747-749.
 
\noindent
[S2] D. Shanks, {\it Two theorems of Gauss}, Pacific. J. Math., 
{\bf 8}(1958), 609-612.
 
\noindent
[WZ] H.S. Wilf and D. Zeilberger, {\it An algorithmic proof theory
for multisum/integral (ordinary and "q")
hypergeometric identities}, Invent. Math. {\bf 108} (1992),
575-633.
 
\noindent
[Z] D. Zeilberger, {\it The method of creative telescoping for q-series},
in preparation.
 
\medskip
\noindent
Feb. 1992 ; Revised: Aug. 1992.
\bye